\documentclass[12pt,reqno]{amsart}
\usepackage[utf8]{inputenc}
\usepackage{amsmath,amsthm,amssymb}
\usepackage{mathtext}
\usepackage[T1,T2A]{fontenc}
\usepackage[english]{babel}
\usepackage{tikz}

\newcommand{\KK}{\mathbb K}

\renewcommand{\phi}{\varphi}

\usepackage{graphicx}
\usepackage[dvipsnames]{xcolor}
\usepackage[colorlinks=true, allcolors=red]{hyperref}
\usepackage{cite}

\theoremstyle{plain}
\newtheorem{theorem}{Theorem}

\theoremstyle{remark}
\newtheorem{remark}{Remark}
\theoremstyle{definition}

\newtheorem{example}{Example}
\newtheorem{question}{Question}

\title{}

\author{Veronika Kikteva}
\address{HSE University, Faculty of Computer Science, 11 Pokrovsky Bulvar, Moscow, 109028, Russia}
\email{VVKikteva@yandex.ru}

\thanks{The work was supported by the Foundation for the Advancement of Theoretical Physics and Mathematics “BASIS”}



\begin{document}
\title{Rigidity conditions for binomials}
\maketitle
\begin{abstract}
We consider the following three conditions for polynomials in several variables. Condition \textbf{(1)} holds for a polynomial if, for every affine integral domain, the following is true: whenever the result of substituting algebraically independent elements into the polynomial lies in the kernel of a locally nilpotent derivation, the elements themselves must lie in that kernel. A polynomial satisfies condition \textbf{(2)} if it does not belong to the kernel of any nonzero LND of the polynomial algebra. It satisfies condition \textbf{(3)} if the quotient algebra modulo the ideal generated by the polynomial is rigid, that is, it admits no nontrivial LNDs. For irreducible binomials, we prove that all three conditions are equivalent. For reducible binomials, we establish all existing implications between these conditions. 
\end{abstract}

\section{Introduction}
Let $\mathbb{K}$ be an algebraically closed field of characteristic zero. We denote by $\mathbb G_a$ the additive group of $\mathbb{K}$. All rings considered in this paper are associative, commutative, and have an identity element. Let $A$ be an algebra over the field $\mathbb K$. A \textit{derivation} of $A$ is a linear map $\delta$: $A\to A$ satisfying the Leibniz rule, that is, for any elements $a,b\in A$, we have $\delta(ab)=a\delta(b)+b\delta(a)$. A derivation $\delta$ of $A$ is called \textit{locally nilpotent} (\textit{LND} for short) if for every~$a\in A$ there exists a natural number $m$ such that $\delta ^m (a)=0$. We denote the set of locally nilpotent derivations of the algebra~$A$ by $\mathrm{LND}(A)$.

Given a locally nilpotent derivation $\delta$ on $A$, one can construct an algebraic $\mathbb{G}_a$-action on~$A$ by setting
$$g\cdot a=\mathrm{exp}(g\delta)(a)=\sum_{i=0}^{\infty}\frac{g^i}{i!}\delta^i(a)$$
for all $g\in \mathbb{K}$ and $a\in A$. The above correspondence is one-to-one. Moreover, the kernel of a locally nilpotent derivation coincides with the set of invariants of the corresponding $\mathbb{G}_a$-action.

We say that a polynomial $P\in\mathbb{K}[X_1,\dots,X_m]$ satisfies \textbf{(1)} if for any affine $\mathbb{K}$-domain $B$, any algebraically independent elements $a_1,\dots,a_m\in B$, and any $\delta\in \mathrm{LND}(B)$, the condition $P(a_1,\dots,a_m)\in \mathrm{Ker}(\delta)$ implies
${a_1,\dots,a_m\in \mathrm{Ker}(\delta)}$.

We give several examples illustrating condition~\textbf{(1)}. The kernel of an LND is factorially closed, that is, for any $\mathbb K$-domain $B$, any LND $\delta$ on $B$, and any nonzero elements $a,b\in B$ the condition $ab\in \mathrm{Ker}(\delta)$ implies $a,b\in\mathrm{Ker}(\delta)$; see~\cite[Principle 1 (a)]{F}. Hence, the monomials $X_1^{\alpha_1}\dots X_m^{\alpha_m}$, where $\alpha_i\geq 1$ for each $i$, satisfy condition \textbf{(1)}. Another class of polynomials satisfying condition \textbf{(1)} consists of polynomials of the form ${c_1X_1^{\alpha_1}+c_2X_2^{\alpha_2}}$, where ${\alpha_1,\alpha_2\geq 2}$, ${c_1,c_2\in\mathbb{K}}$, and $m=2$, as follows from~\cite[Theorem 2.50]{F}. A generalization of this result to $k$-nomials of the form
$$c_1 X_{11}^{\alpha_{11}}\dots X_{1n_1}^{\alpha_{1n_1}}+\dots +c_k X_{k1}^{\alpha_{k1}}\dots X_{kn_k}^{\alpha_{kn_k}},$$
where $k\geq 3$ and $\sum 1/\alpha_{ij}\leq 1/(k-1)$, in the case $m=n_1+\dots +n_k$, was obtained in~\cite[Theorem 3.3]{K}. Examples of polynomials that do not satisfy condition \textbf{(1)} in the case ${m>1}$ include coordinate functions, as well as polynomials that can be included in some coordinate system.

We say that a polynomial $P\in \mathbb{K}[X_1,\dots,X_m]$ satisfies condition \textbf{(2)} if it does not belong to the kernel of any nonzero locally nilpotent derivation $D\in\mathrm{LND}(\mathbb{K}[X_1,\dots, X_m])$.

\begin{remark}
Condition \textbf{(2)} follows from condition \textbf{(1)}. Indeed, let $P\in\mathbb{K}[X_1,\dots,X_m]$ satisfy condition~\textbf{(1)}. Take $B=\mathbb{K}[X_1,\dots,X_m]$. The elements $a_i:=X_i$ are algebraically independent. If $P$ belongs to the kernel of ${D\in\mathrm{LND}(\mathbb{K}[X_1,\dots,X_m])}$, then condition~\textbf{(1)} implies that $X_i\in\mathrm{Ker}(D)$, and hence $D$ is the zero derivation.
\end{remark}

The following question was posed in the book~\cite{F}.
\begin{question}\cite[Question 11.9]{F}
Does condition \textbf{(1)} follow from condition \textbf{(2)}?
\end{question}

In the case $m=2$ the answer is positive for irreducible polynomials. In~\cite[Corollary 2.1]{KML}, it was proved that any irreducible polynomial $P\in\mathbb{C}[x,y]$ that is not a coordinate satisfies condition \textbf{(1)}. For an irreducible polynomial in the case $m=2$, the condition of not being a coordinate is equivalent to condition~\textbf{(2)}; see~\cite[Section 11.7]{F}.

The present paper proves the implication \textbf{(2)}$\implies$\textbf{(1)} for binomials.

\begin{theorem}
\label{thmbinom21}
Let $P\in\mathbb{K}[X_1,\dots,X_m]$ be of the form
$$P=X_1^{\alpha_1}\dots X_m^{\alpha_m}+cX_1^{\beta_1}\dots X_m^{\beta_m},$$
where $\alpha_i,\beta_j\geq 0$ and $c\in\mathbb{K}^\times$. Suppose that $P$ satisfies condition~\textbf{(2)}. Then $P$ satisfies condition \textbf{(1)}.
\end{theorem}

An algebra is called \textit{rigid} if it admits no nonzero locally nilpotent derivations. An affine algebraic variety $X$ is called \textit{rigid} if its coordinate ring is rigid. This condition is equivalent to the absence of subgroups isomorphic to~$\mathbb G_a$ in the automorphism group of $X$. There are criteria for the rigidity of various classes of varieties; see, for example,~\cite[Theorem 2, Remark 3]{BG} for toric varieties,~\cite[Theorem 1]{A},~\cite[Theorem 3]{Ga}, and~\cite[Theorem 3]{EGS} for trinomial varieties, and~\cite[Main Theorem, Corollary]{CD} for Pham--Brieskorn varieties for $n=4$.

In general, the automorphism groups of affine algebraic varieties are not algebraic. For example, if a variety admits a nontrivial $\mathbb G_a$-action and has dimension at least 2, then even the identity component of its automorphism group in the sense of~\cite{R} is not an algebraic variety.

For rigid affine algebraic varieties, there are results describing the structure of their automorphism groups. It has been conjectured that the identity component of the automorphism group of a rigid affine variety is an algebraic torus; see~\cite[Conjecture 1.0.1]{PZ}. This conjecture has been proved for normal surfaces; see~\cite[Theorem 1.0.3]{PZ}. It has also been established for not necessarily normal toric varieties and for rational normal affine varieties without nonconstant invertible regular functions, with a finitely generated divisor class group, and admitting a torus action of complexity 1; see~\cite[Remark 6.1, Theorem 6.4]{BorG}. It was also proved in~\cite[Theorem 2.1]{AG} that the automorphism group of a rigid variety contains a unique maximal torus.

For some classes of rigid varieties, their automorphism groups have been described; see, for example,~\cite[Theorem 3]{BG} for toric varieties,~\cite[Proposition 2, Theorem 3]{AG} for trinomial hypersurfaces, and~\cite[Theorem 2.7]{T} for hypersurfaces defined by relations in which each variable occurs in a unique monomial.

We say that a polynomial $P\in\mathbb{K}[X_1,\dots,X_m]$ satisfies condition~\textbf{(3)} if the quotient algebra $\mathbb{K}[X_1,\dots,X_m]/(P)$ is rigid.

\begin{remark}
Condition \textbf{(3)} implies condition \textbf{(2)}. Suppose that a polynomial $P$ belongs to the kernel of a nonzero locally nilpotent derivation $D\in \mathrm{LND}(\mathbb{K}[X_1,\dots,X_m])$. Choose the maximal $k\in\mathbb{N}\cup {0}$ such that $D(X_1),\dots,D(X_m)\in(P^k)$, and replace $D$ by the derivation~$\widetilde{D}$ acting on the generators according to the rule $X_i\mapsto D(X_i)/P^k$.

The property of being locally nilpotent and nontrivial, as well as the condition ${P\in\mathrm{Ker}(D)}$, is preserved; see~\cite[Principle 7]{F}. Consider the algebra
$$A:=\mathbb{K}[X_1,\dots,X_m]/(P)=\mathbb{K}[x_1,\dots,x_m],$$
and define a derivation on $A$ by its action on the generators:
$$\delta: x_i\mapsto \widetilde{D}(X_i)+(P).$$
The resulting derivation is nonzero and locally nilpotent. Hence, the algebra $A$ is not rigid.
\end{remark}

The converse implication is proved in this paper for irreducible binomials.

\begin{theorem}
\label{thmbinom23}
For irreducible binomials, condition \textbf{(2)} implies condition \textbf{(3)}.
\end{theorem}

It is worth noting that, for irreducible binomials $P$, the affine variety $\mathbb{V}(P)$ is toric; see~\cite[Proposition 1.1.11, Theorem 1.1.17]{CLS}. A criterion for the rigidity of affine toric varieties was obtained in~\cite{BG} and formulated in terms of the combinatorial data associated with the toric variety. However, we will not interpret condition \textbf{(2)} in terms of toric data; the present paper does not rely on the results of~\cite{BG}. More detailed information on binomial ideals can be found in~\cite{ES, HHO}.

For reducible binomials, the implication \textbf{(2)}$\implies$\textbf{(3)} may fail, as shown by the following example.

\begin{example}
The polynomial $X(XY-1)\in\mathbb{K}[X,Y]$ does not belong to the kernel of any nonzero LND on $\mathbb K[X,Y]$. Indeed, if $X(XY-1)\in \mathrm{Ker} (D)$, then, since the kernel of an LND is factorially closed, we have $X,XY-1\in\mathrm{Ker} (D)$, and hence $X,Y\in \mathrm{Ker} (D)$. However, the corresponding quotient algebra
$$\mathbb{K}[X,Y]/(X(XY-1))=\mathbb{K}[x,y]$$
is not rigid: it admits a derivation $\delta$ defined on the generators by
$$\delta(x)=0,\ \delta(y)=xy-1.$$
This derivation is well-defined, since
$$\delta : x(xy-1)\mapsto x^2(xy-1)=0,$$
and it is nonzero and locally nilpotent, since
$$\delta: x\mapsto 0,\ y\mapsto xy-1\mapsto x(xy-1)=0.$$
Thus, the polynomial $X(XY-1)\in\mathbb{K}[X,Y]$ satisfies condition \textbf{(2)} but does not satisfy condition \textbf{(3)}.
\end{example}

The author is grateful to Sergey Gaifullin for constant
attention to this work.

\section{Preliminaries}

\label{sectprelim}

We recall the basic definitions related to locally nilpotent derivations. Detailed information and proofs can be found in the book~\cite{F}.

Let $B$ be a $\mathbb K$-domain. Suppose that $B$ is graded by an abelian group $G$, that is, there is a family $\{B_g\}_{g\in G}$ of subspaces of $B$ such that
$$B=\bigoplus_{g\in G}B_g,\text{ and } B_{g_1}B_{g_2}\subseteq B_{g_1+g_2} \text{ for all }g_1,g_2\in G .$$
A locally nilpotent derivation $\delta$ on a $G$-graded algebra $B$ is called $G$-\textit{homogeneous} if there exists an element $d\in G$ such that $\delta (B_g)\subseteq B_{g+d}$ for all $g\in G$. If $\delta \neq 0$, then $d$ is uniquely determined, is called the \textit{degree} of $\delta$, and is denoted by $\mathrm{deg}(\delta)$.

Every nonzero LND $\delta$ on a $G$-graded algebra $B$ can be decomposed into a sum of nonzero $G$-homogeneous derivations $\delta = \sum_{i\in I} \delta_i$ for a finite subset $I\subseteq G$, where $\mathrm{deg}(\delta_i)=i$. If $G=\mathbb{Z}$, the components of minimal and maximal degree are also locally nilpotent, that is, $\delta_m,\delta_n\in\mathrm{LND}(B)$, where $m=\mathrm{min}\,I,\ n=\mathrm{max}\,I$; see~\cite[Corollary 1.31]{F} and~\cite[Principle II]{E}. Thus, a $\mathbb{Z}$-graded algebra admits a nonzero LND if and only if it admits a nonzero $\mathbb{Z}$-homogeneous LND. Applying this result $n$ times, we obtain that a $\mathbb{Z}^n$-graded algebra admits a nonzero LND if and only if it admits a nonzero $\mathbb{Z}^n$-homogeneous LND; see also~\cite[Lemma 2.4]{G}.

A \textit{degree function} on $B$ is a map $\mathrm{deg}:B\to G\cup \{-\infty\}$, where $G$ is a linearly ordered abelian group, satisfying the following conditions for all $a,b\in B$:
\begin{enumerate}
\item $\mathrm{deg}(a)=-\infty$ if and only if $a=0$;
\item $\mathrm{deg}(ab)=\mathrm{deg}(a)+\mathrm{deg}(b)$;
\item $\mathrm{deg}(a+b)\leq \mathrm{max}\{\mathrm{deg}(a), \mathrm{deg}(b)\}$.
\end{enumerate}
Let $\delta$ be a locally nilpotent derivation on $B$. Then $\delta$ defines a degree function on $B$ with values in $\mathbb N\cup \{0, -\infty\}\subseteq\mathbb Z\cup \{ -\infty\}$. For every nonzero $b\in B$, set
$$\mathrm{deg}_{\delta}(b)=\mathrm{min}\{n\in \mathbb N\cup \{0\} \mid \delta ^{n+1}(b)=0\}$$
and set $\mathrm{deg}_{\delta}(0)=-\infty$.

An element $t$ of an algebra $B$ is called a \textit{local slice} for $\delta\in\mathrm{LND}(B)$ if $\delta^2(t)=0$, but $\delta(t)\neq 0$. Every nonzero LND admits a local slice. For each $f\in B\backslash\{0\}$, let $B_f$ denote the localization of $B$ with respect to the multiplicatively closed set $\{f^i\}_{i\geq 0}$.

Let $A:=\mathrm{Ker}(\delta)$ and let $t\in B$ be a local slice for $\delta\in\mathrm{LND}(B)$. The local slice theorem states that $B_{\delta (t)}=A_{\delta (t)}[t]$; see~\cite[Principle 11 (d)]{F}. Consequently, $B$ can be embedded into a polynomial ring in $t$ over the field of fractions of the kernel of the LND. In this setting, $\delta$ coincides with the derivation with respect to the variable~$t$, and the degree of elements of~$B$ as polynomials in $t$ coincides with the degree induced by $\delta$: ${\mathrm{deg}_t=\mathrm{deg}_\delta}$. In particular, an element of $B$ belongs to the kernel of $\delta$ if and only if it is constant as a polynomial in~$t$.

Now let $F\in\mathbb{L}[T]$ be a polynomial in one variable $T$ over an algebraically closed field $\mathbb L$ of characteristic zero. Denote by $N(F)$ the number of distinct roots of $F$ in $\mathbb{L}$. Further we will use the Mason--Stothers theorem.

\begin{theorem}\cite{S}
\label{thmMS}
Let $A,B,C\in \mathbb L[T]$ be pairwise coprime polynomials in $T$ over an algebraically closed field of characteristic zero, with at least one of them nonconstant. If $A+B+C=0$, then
$$\max\{ \mathrm{deg}_T(A), \mathrm{deg}_T(B), \mathrm{deg}_T(C)\}\leq N (ABC)-1.$$
\end{theorem}

\section{Proof of Theorem~\ref{thmbinom21}}

In this section, we prove that if a polynomial $P\in\mathbb{K}[X_1,\dots,X_m]$ of the form
$$P=X_1^{\alpha_1}\dots X_m^{\alpha_m}+cX_1^{\beta_1}\dots X_m^{\beta_m},$$
where $\alpha_i,\beta_j\geq 0$ and $c\in\mathbb{K}^\times$, does not belong to the kernel of any nonzero LND of the algebra $\mathbb{K}[X_1,\dots, X_m]$, then $P$ satisfies condition \textbf{(1)}.

The polynomial $P$ depends non-trivially on each of the variables $X_1,\dots,X_m$, that is, $\alpha_i+\beta_i>0$ for every $i=1,\dots,m$, since otherwise $P$ belongs to the kernel of the LND~$\partial/\partial X_i$.

If $\sum \beta_i=0$, then $P=X_1^{\alpha_1}\dots X_m^{\alpha_m}+c$. In this case, $\alpha_i>0$ for every $i=1,\dots,m$. Let~$B$ be any affine $\mathbb{K}$-domain, let $a_1,\dots,a_m\in B$ be algebraically independent elements, and let $\delta\in \mathrm{LND}(B)$. If
$a_1^{\alpha_1}\dots a_m^{\alpha_m}+c\in \mathrm{Ker}(\delta)$, then $a_1^{\alpha_1}\dots a_m^{\alpha_m}\in \mathrm{Ker}(\delta)$, since $\mathbb{K}\subseteq \mathrm{Ker}(\delta)$. Hence $a_1,\dots,a_m\in \mathrm{Ker}(\delta)$, since the kernel of an LND is factorially closed; see \cite[Principle 1(a)]{F}.

Thus, we may assume that $\sum \alpha_i,\sum\beta_j\neq 0$. Observe that
$$P=X_1^{\widetilde{\alpha_1}} \dots X_m^{\widetilde{\alpha_m}} (X_1^{\widetilde{\beta_1}} \dots X_m^{\widetilde{\beta_m}}+cX_1^{\widetilde{\gamma_1}}\dots X_m^{\widetilde{\gamma_m}}),$$
where $\widetilde{\alpha_i}=\min(\alpha_i,\beta_i)$, $\widetilde{\beta_i}=\alpha_i-\widetilde{\alpha_i}$, $\widetilde{\gamma_i}=\beta_i-\widetilde{\alpha_i}$, and, for every $i=1,\dots,m$, at least one of the three numbers $\widetilde{\alpha_i},\widetilde{\beta_i},\widetilde{\gamma_i}$ is nonzero, while one of $\widetilde{\beta_i},\widetilde{\gamma_i}$ is zero.

After relabeling the variables, we divide $X_1,\dots,X_m$ into three disjoint groups. Let $X_1,\dots,X_k$ be the variables for which $\widetilde{\alpha_i}\neq 0$, let $X_i=:Y_i$ if ${\widetilde{\alpha_i}=0},\  \widetilde{\beta_i}\neq 0$, and let $X_i=:Z_i$ if $\widetilde{\alpha_i}=0,\ \widetilde{\gamma_i}\neq 0$. Denote the corresponding exponents by ${\xi_i},{\eta_i},{\zeta_i}$, respectively. Note that ${\xi_i},{\eta_i},{\zeta_i}>0$. Then $P$ has the form
$$P=X_1^{{\xi_1}}\dots X_k^{{\xi_k}} ( P_1 Y_1^{{\eta_1}}\dots Y_r^{{\eta_r}}+P_2Z_1^{{\zeta_1}}\dots Z_s^{{\zeta_s}}),$$
where $P_1,P_2$ are monomials in $X_1,\dots,X_k$.

We show that ${\eta_i},{\zeta_i}\geq 2$. If ${\eta_1}=1$, then $P$ belongs to the kernel of the locally nilpotent derivation $D$ defined on $Y_1,Z_1$ by
$$D:\ Y_1 \mapsto P_2 {\zeta_1} Z_1^{{\zeta_1}-1} Z_2^{{\zeta_2}}\dots Z_s^{{\zeta_s}},$$
$$D:\ Z_1\mapsto -P_1  Y_2^{{\eta_2}}\dots Y_r^{{\eta_r}} ,$$
and equal to zero on all other variables. Indeed,
$$D(P)=X_1^{{\xi_1}}\dots X_k^{{\xi_k}} ( P_1 D(Y_1)Y_2^{{\eta_2}}\dots Y_r^{{\eta_r}}+P_2D(Z_1^{{\zeta_1}})Z_2^{{\zeta_2}}\dots Z_s^{{\zeta_s}})=$$
$$=X_1^{{\xi_1}}\dots X_k^{{\xi_k}} ( P_1P_2 
{\zeta_1}
 Y_2^{{\eta_2}}\dots Y_r^{{\eta_r}}
 Z_1^{{\zeta_1}-1} Z_2^{{\zeta_2}} \dots Z_s^{{\zeta_s}}
-$$
$$-P_1P_2
{\zeta_1}
Y_2^{{\eta_2}}\dots Y_r^{{\eta_r}} Z_1^{{\zeta_1}-1}
Z_2^{{\zeta_2}}\dots Z_s^{{\zeta_s}})=0.$$

Now let $B$ be an affine $\mathbb{K}$-domain, let $a_1,\dots,a_k,b_1,\dots,b_r,c_1,\dots,c_s\in B$ be algebraically independent elements, and let $\delta\in \mathrm{LND}(B)$. Suppose that
$$P(a_1,\dots,a_k,b_1,\dots,b_r,c_1,\dots,c_s)\in \mathrm{Ker}(\delta),$$
and we prove that $a_i,b_i,c_i\in\mathrm{Ker}(\delta)$. If $\delta=0$, the statement is obvious. Hence, we may assume that $\delta\in \mathrm{LND}(B)\backslash\{0\}$. The elements $a_1,\dots,a_k$ belong to the kernel of $\delta$, since the kernel of an LND is factorially closed. Consequently, the monomials $P_1$ and $P_2$ in these elements also belong to the kernel of $\delta$. For convenience, below we write $p$ and $p_i$ for $P(a_1,\dots,a_k,b_1,\dots,b_r,c_1,\dots,c_s)$ and $P_i(a_1,\dots,a_k)$, respectively. Thus,
$$p=a_1^{{\xi_1}}\dots a_k^{{\xi_k}} ( p_1 b_1^{{\eta_1}}\dots b_r^{{\eta_r}}+p_2c_1^{{\zeta_1}}\dots c_s^{{\zeta_s}})\in \mathrm{Ker}(\delta).$$
By factorial closedness of the kernel of an LND, it follows that
$$p_1 b_1^{{\eta_1}}\dots b_r^{{\eta_r}}+p_2c_1^{{\zeta_1}}\dots c_s^{{\zeta_s}}\in \mathrm{Ker}(\delta),$$
where ${\eta_i},{\zeta_i}\geq 2$.

Since $\delta$ is a nonzero LND, it admits a local slice $t$. Embed $B$ into $\mathbb{L}[t]$, where $\mathbb{L}=\mathrm{Quot}(\mathrm{Ker}(\delta))$, so that ${\mathrm{deg}_t=\mathrm{deg}_\delta=:\mathrm{deg}}$. Thus,
$$p_1 b_1^{{\eta_1}}\dots b_r^{{\eta_r}}+p_2c_1^{{\zeta_1}}\dots c_s^{{\zeta_s}}=:c'\in \mathbb{L}^\times.$$
The left-hand side cannot be zero, since the elements $a_1,\dots,a_k$, $b_1,\dots,b_r$, $c_1,\dots,c_s$ are algebraically independent over $\mathbb K$. Suppose that not all elements $b_i,c_j$ are constant as polynomials in $t$. Then we apply the Mason--Stothers theorem to the elements $p_1 b_1^{{\eta_1}}\dots b_r^{{\eta_r}}$, $p_2c_1^{{\zeta_1}}\dots c_s^{{\zeta_s}}$, and $-c'$. The theorem can be applied, since these polynomials are pairwise coprime. Indeed, if the first two had a nontrivial common divisor, then this divisor would divide the constant $c'$.

Let $\overline{\mathbb{L}}$ denote the algebraic closure of the field $\mathbb L$, and let $N(f)$ denote the number of distinct roots of a polynomial $f\in\overline{\mathbb{L}}[t]$ in $\overline{\mathbb{L}}$. We have
$$\max\{ \mathrm{deg}(p_1 b_1^{{\eta_1}}\dots b_r^{{\eta_r}}), \mathrm{deg}(p_2c_1^{{\zeta_1}}\dots c_s^{{\zeta_s}}),\mathrm{deg}(-c')\}\leq$$
$$\leq  N (p_1 b_1^{{\eta_1}}\dots b_r^{{\eta_r}}p_2c_1^{{\zeta_1}}\dots c_s^{{\zeta_s}}(-c'))-1,$$
that is,
$$\max\{ \mathrm{deg}(b_1^{{\eta_1}}\dots b_r^{{\eta_r}}), \mathrm{deg}(c_1^{{\zeta_1}}\dots c_s^{{\zeta_s}})\}\leq N (b_1\dots b_r c_1\dots c_s)-1.$$
Since ${\eta_i},{\zeta_j}\geq 2$ for all $i,j$, we have
$$\max\{ 2 \mathrm{deg}(b_1)+\dots +2\mathrm{deg}(b_r),2\mathrm{deg}(c_1)+\dots + 2 \mathrm{deg}(c_s)\}  \leq $$
$$ \leq \max\{ \mathrm{deg}(b_1^{{\eta_1}}\dots b_r^{{\eta_r}}), \mathrm{deg}(c_1^{{\zeta_1}}\dots c_s^{{\zeta_s}})\} \leq N (b_1\dots b_r c_1\dots c_s)-1 \leq $$
$$ \leq \mathrm{deg}(b_1)+\dots +\mathrm{deg}(b_r) + \mathrm{deg}(c_1)+\dots +\mathrm{deg}(c_s)-1.$$
Hence,
$$2\max \{\mathrm{deg}(b_1)+\dots +\mathrm{deg}(b_r), \mathrm{deg}(c_1)+\dots +\mathrm{deg}(c_s)\} \leq $$
$$\leq \mathrm{deg}(b_1)+\dots +\mathrm{deg}(b_r) +\mathrm{deg}(c_1)+\dots +\mathrm{deg}(c_s)-1,$$
which is a contradiction. Therefore, the elements $b_i,c_j$ are constant as polynomials in $t$ and belong to the kernel of $\delta$.

\section{Proof of Theorem~\ref{thmbinom23}}

In this section, we prove that if the irreducible polynomial
$$P=X_1^{\alpha_1}\dots X_r^{\alpha_r}+cY_1^{\beta_1}\dots Y_s^{\beta_s},$$
where $\alpha_i,\beta_j\geq 0$ and $c\in\mathbb{K}^\times$, does not belong to the kernel of any nonzero LND of the polynomial algebra~$\mathbb K[X_1,\dots,X_r,Y_1,\dots, Y_s]$, then the quotient algebra by the ideal generated by $P$ is rigid.

We may assume that either
$$P=X_1^{\alpha_1}\dots X_m^{\alpha_m}+c\in\mathbb K[X_1,\dots,X_m],$$
where $c\in \mathbb{K}^{\times}$ and all $\alpha_i>0$, or
$$P=X_1^{\alpha_1}\dots X_r^{\alpha_r}+Y_1^{\beta_1}\dots Y_s^{\beta_s}\in\mathbb K[X_1,\dots,X_r,Y_1,\dots,Y_s],$$
where $\alpha_i,\beta_j\geq 2$, since otherwise one can construct a locally nilpotent derivation of the polynomial algebra that vanishes on $P$, as in the proof of Theorem~\ref{thmbinom21}.

\subsection{The case $P=X_1^{\alpha_1}\dots X_m^{\alpha_m}+c$}

Consider the domain
$$B=\mathbb{K}[X_1,\dots,X_m]/(P)=\mathbb{K}[x_1,\dots,x_m]$$
and an arbitrary locally nilpotent derivation $\delta$ on $B$. We have
$$x_1^{\alpha_1}\dots x_m^{\alpha_m}=-c\in \mathrm{Ker}(\delta) \backslash\{0\},$$
and hence $x_1,\dots,x_m\in \mathrm{Ker}(\delta)$, since the kernel of an LND is factorially closed. Therefore, $\delta=0$, and the algebra $B$ is rigid.

\subsection{The case $P=X_1^{\alpha_1}\dots X_r^{\alpha_r}+Y_1^{\beta_1}\dots Y_s^{\beta_s}$}

We prove that the algebra
$$B=\mathbb{K}[X_1,\dots,X_r,Y_1,\dots,Y_s]/(P)=\mathbb{K}[x_1,\dots,x_r,y_1,\dots,y_s]$$
is rigid.

Set $q=\binom{r}{2}+\binom{s}{2}$ and label the first $\binom{r}{2}$ coordinates in $\mathbb{Z}^q$ by the two-element subsets $\{i,j\}$ for $\{i,j\}\subseteq \{1,\dots,r\}$ and the remaining coordinates by the two-element subsets $\{k,l\}\subseteq \{1,\dots,s\}$. Let the component of $\deg(x_i)\in \mathbb Z^q$ corresponding to $\{i,j\}$ be $\alpha_j$ if $i<j$, and $-\alpha_j$ otherwise, with all other components equal to~$0$. Similarly, at the position $\{k,l\}$, let $\deg(y_k)$ be $\beta_l$ if $k<l$ and $-\beta_l$ if $k>l$, with all other components equal to~$0$. It is easy to check that
$$\deg(x_1^{\alpha_1}\dots x_r^{\alpha_r})=\deg(y_1^{\beta_1}\dots y_s^{\beta_s})=0.$$

If the algebra $B$ is not rigid, then it admits a nonzero $\mathbb{Z}^q\text{-homogeneous}$ LND~$\delta$. Let its degree be
$$\deg(\delta)=(d_1,\dots,d_q)=(d_{\{i,j\}},d_{\{ k,l\}})_{\{i,j\}\subseteq \{1,\dots,r\},\{k,l\}\subseteq \{1,\dots,s\}}.$$

For each two-element subset $\{i,j\}\subseteq \{1,\dots,r\}$ with $i<j$, consider the projection homomorphism $\pi_{\{i,j\}}:\mathbb{Z}^q\to \mathbb{Z}$ onto the $\{i,j\}$-coordinate. Then $B$ is also a $\mathbb{Z}$-graded algebra, and $\delta$ is a $\mathbb{Z}$-homogeneous LND; see~\cite[Lemma 2.5]{G}. We denote the degree corresponding to this grading by $\deg_{\{i,j\}}$. For the generators of $B$, we have
$$\deg_{\{i,j\}}(x_i)=\alpha_j>0,\ \deg_{\{i,j\}}(x_j)=-\alpha_i<0,$$
while the degrees of the other coordinate functions are zero, and $\deg_{\{i,j\}}(\delta)=d_{\{i,j\}}$. If $d_{\{i,j\}}\geq 0$, then
$$\deg_{\{i,j\}}(\delta(x_i)) = \alpha_j+d_{\{i,j\}}>0.$$
Hence, $x_i \mid \delta(x_i)$, and by~\cite[Corollary 1.23]{F} we have $x_i\in\mathrm{Ker}(\delta)$. Similarly, if $d_{\{i,j\}}\leq 0$, then
$$\deg_{\{i,j\}}(\delta(x_j)) = -\alpha_i+d_{\{i,j\}}<0$$
and $x_j\in \mathrm{Ker}(\delta)$.

Thus, for any distinct $i,j\in\{1,\dots,r\}$, at least one of the elements $x_i$ and $x_j$ belongs to the kernel of~$\delta$. Hence, at most one of $x_1,\dots,x_r$ does not belong to $\mathrm{Ker}(\delta)$. The analogous statement holds for $y_1,\dots, y_s$. Without loss of generality, we may assume that all the listed elements except $x_1,y_1$ belong to the kernel of~$\delta$.

Extend the field $\KK$ by adjoining $X_2,\dots,X_r,Y_2,\dots,Y_s$, and consider its algebraic closure
$$\mathbb{L}:= \overline{\mathbb{K}(X_2,\dots,X_r,Y_2,\dots,Y_s)}.$$
The algebra $B$ can be embedded into
$$B'=\mathbb{L}[x_1,y_1]=\mathbb{L}[X_1,Y_1]/(P_1 X_1^{\alpha_1}+P_2 Y_1^{\beta_1}),$$
where
$$P_1:=X_2^{\alpha_2}\dots X_r^{\alpha_r},\ P_2:=Y_2^{\beta_2}\dots Y_s^{\beta_s}\in \mathbb{L},$$
and $\delta$ extends to a nontrivial LND on $B'$. Since $\alpha_1,\beta_1\geq2$ by the argument above, Lemma~4.5 of~\cite{G} implies that $B'$ is rigid.

Thus, every $\mathbb{Z}^q$-homogeneous LND on $B$ is zero, and hence the algebra $B$ is rigid.

\end{document}